\documentclass{elsart}

\usepackage{amssymb,amsmath}

\newtheorem{theorem}{Theorem}[section]

\newtheorem{proposition}[theorem]{Proposition}
\newtheorem{corollary}[theorem]{Corollary}

\newtheorem{question}[theorem]{Question}
\newtheorem{example}[theorem]{Example}

\journal{}
\begin{document}

\begin{frontmatter}

\title{Diagonals of separately absolutely continuous mappings coincide with the sums of absolute convergent series of continuous functions}

\author{Olena Karlova, Volodymyr Mykhaylyuk, Oleksandr Sobchuk}

\address{Chernivtsi National University, Department of Mathematical Analysis,
Kotsjubyns'koho 2, Chernivtsi 58012, Ukraine\\maslenizza.ua@gmail.com}

\begin{abstract}
We prove that, for an interval $X\subseteq \mathbb R$ and a normed space $Z$ diagonals of separately absolute continuous mappings $f:X^2\to Z$ are exactly mappings \mbox{$g:X\to Z$} which are the sums of absolute convergent series of continuous functions.
\end{abstract}

\end{frontmatter}

\maketitle

\section{Introduction}

Let $f:X^2\to Z$ be a mapping. We call a mapping  $g:X\to Z$, $g(x)=f(x,x)$, {\it the diagonal of $f$.}

Investigations of diagonals of separately continuous functions $f:X^2\to\mathbb R$ started in the classical work of R.~Baire \cite{Baire}. He showed that diagonals of separately continuous functions of two real variables are exactly Baire-one functions, i.e. pointwise limits of continuous functions. His result was generalized by A.~Lebesgue and H.~Hahn for real-valued functions of several real variables (see \cite{Leb1,Leb2,Ha}).

Since the second half the XX-th century, Baire classification of separately continuous mappings and their {analogues} is intensively studied by many mathe\-ma\-ticians (see \cite{Mo,Ru,Ve,MMMS,B,Bur}). The inverse problem on construction of separately continuous functions  with {a} given diagonal was solved in \cite{MMS1}. In \cite{My} it was shown that for any topological space $X$ and a function $g:X\to\mathbb R$ of the $(n-1)$-th Baire class there exist a separately continuous function  $f:X^n\to\mathbb R$ with the diagonal $g$.

In \cite{MMP} the diagonal variant of the problem of Eidelheit from {the} famous ``Scottish Book'' on a composition of absolutely continuous functions was investigated. It was proved that there exists a separately absolutely continuous function $f:[0,1]^2\to\mathbb R$ such that its partial derivatives  $f'_x$ and $f'_y$ in {the degree} $p$ are integrable on $[0,1]^2$ for every $p>1$, and such that its diagonal $g$ is not absolutely continuous.

{The following question naturally arises}.

\begin{question}\label{q:1.1}
Find necessary and sufficient conditions on a function $g: [0,1] \to \mathbb R$
 under  which  there  is  a  separately  absolutely continuous function
 $f:[0,1]^2\to \mathbb R$ with the diagonal $g$.
\end{question}

In the given paper we prove that for any interval $X\subseteq \mathbb R$  and a normed space $Z$ diagonals of separately absolute continuous mappings $f:X^2\to Z$ are exactly mappings \mbox{$g:X\to Z$} which are the sums of absolute convergent series of continuous functions.

\section{Preliminaries}

For topological spaces $X$, $Y$ and $Z$ a mapping $f:X\times Y\to Z$ which is continuous with respect to every variable is called  {\it separately continuous}.

Let $X$ and $Y$ be topological spaces. A mapping $f:X\to Y$ is {\it a mapping of the first Baire class} or {\it a Baire-one mapping} if there exists a sequence $(f_n)^{\infty}_{n=1}$ of continuous mappings $f_n:X\to Y$ which {pointwise converges to} $f$ on $X$.

For a metric space $X$ by $|\cdot - \cdot|_X$ we denote {the} metric on this space.

Let $X\subseteq \mathbb R$ be an interval and let  $Z$ be a metric space. A mapping $f:X\to Z$ is called {\it absolutely continuous} if for an arbitrary $\varepsilon>0$ there exists $\delta>0$ such that for every collection $a_1< b_1\leq a_2<b_2\leq \dots \leq a_n<b_n$ of elements $a_1, b_1, \dots ,a_n, b_n\in X$ with $\sum\limits_{k=1}^{n}(b_k-a_k)<\delta$ the inequality $\sum\limits_{k=1}^{n}|f(b_k)-f(a_k)|_Z<\varepsilon$ holds. Let, moreover, $Y$ be an interval. A mapping $f:X\times Y\to Z$ which is absolutely continuous with respect to each variable is called  {\it separately  absolutely continuous}.

A mapping $f:X\to Z$ has {\it {bounded} variation on an interval $X$} if there exists  $C>0$ such that for any collection $a_1< b_1\leq a_2<b_2\leq \dots \leq a_n<b_n$ of elements $a_1, b_1, \dots ,a_n, b_n\in X$ the inequality $\sum\limits_{k=1}^{n}|f(b_k)-f(a_k)|_Z\leq C$ holds.
Moreover, for an interval $X=[a,b]$ the least upper bound of all values $\sum\limits_{k=1}^{n}|f(b_k)-f(a_k)|_Z$ is called {\it the variation of $f$ on $[a,b]$}.

Let $X\subseteq \mathbb R$ be an interval, let $Z$ be a metric space, let $f:X\to Z$ be a mapping and $x_0\in X$. We say that {\it $f$ {has finite} variation at $x_0$} if there exists a segment $[a,b]\subseteq X$ such that $[a,b]$ is a neighborhood of $x_0$ in $X$ and $f$ {has  finite} variation on $[a,b]$.

Let $X$ and $Z$ be metric spaces and $A\subseteq X$. A mapping $f:X\to Z$ is {\it Lipschitz on a set $A$ with a constant $C\geq 0$} if $|f(x)-f(y)|_Z\leq C|x-y|_X$ for any $x,y\in A$. A mapping $f:X\to Z$ is called {\it Lipschitz on a set  $A$} if there exists $C\geq 0$ such that $f$ is Lipschitz on $A$ with the constant $C$. A mapping $f:X\to Z$ is called {\it Lipschitz (with a constant $C$)} if it is Lipschitz (with a constant $C$) on $X$.

A function $p:X\to\mathbb R$ defined on a vector space $X$ over a field $\mathbb K$ is said to be {\it a pseudo-norm} if for any $x,y\in X$ and $\lambda\in\mathbb K$ with  $|\lambda|\leq 1$ the following conditions hold: $p(x)\geq 0$, moreover, $p(x)=0$ iff $x=0$; $p(\lambda x)\leq p(x)$ and $p(x+y)\leq p(x)+p(y)$.
It is well-known that for any metrizable topological vector space  $X$ there exists a metric $\rho$ on $X$ which generates a topological structure on $X$ and $\rho(x,y)=p(x-y)$ for some pseudo-norm $p$ on $X$.

Let $X$ be a topological space and let $Z$ be a metric space.  A mapping $f:X\to Z$ is said to be {\it an absolute Baire-one mapping} or {\it a mapping of the first absolute Baire class} if there exists a sequence $(f_n)_{n=1}^{\infty}$ of continuous functions $f_n:X\to Z$ such that  $\lim\limits_{n\to\infty}f_n(x)=f(x)$ and \mbox{$\sum\limits_{n=1}^{\infty}|f_{n+1}(x)-f_n(x)|_Z<\infty$} for every $x\in X$. Note that for a normed space $Z$ a mapping $f:X\to Z$ is an absolute Baire-one mapping if and only if $f$ is the sum of an absolute convergent series of continuous functions. It is well-known (see \cite{S}, \cite[\S 41 section 5, \S 42 section 2]{Haus}) that a real valued function is the sum of an absolute convergent series of continuous functions if and only if it is the sum of a lower semicontinuous and an upper semicontinuous functions. Some properties of real-valued absolute Baire-one functions were studied by many mathematicians (see, for example, \cite{Mor},\cite{HOR}).

\section{Necessary conditions on diagonals of separately absolutely continuous mappings}

\begin{proposition}\label{pr:3.1}
  Let $X=\mathbb R$, let $Z$ be a metric topological vector space with the  metric, generated by a pseudo-norm $p$, let $f:X^2\to Z$ be a continuous mapping with respect to the first variable and let $\alpha:X\to (0,+\infty)$ be a continuous function. Then there exist functions $\beta, \gamma:X\to (\frac{1}{2}\alpha(x),\alpha(x))$ and a continuous mapping $g:X\to Z$ such that $$p(g(x)-f(x,x+\beta(x)))\leq p(f(x,x+\beta(x))-f(x,x+\gamma(x)))$$ for every $x\in X$.
\end{proposition}

{\bf Proof.} For every $x\in X$ we write $y_x=x+\frac{3}{4}\alpha(x)$ and, taking into account that  $\alpha$ is continuous at $x$, we choose an open neighborhood $U_x$ of $x$ in $X$ such that $t+\frac{1}{2}\alpha(t)<y_x<t+\alpha(t)$ for every $t\in U_x$.

Let $(V_i:i\in I)$ be an open locally finite refinement of $(U_x:x\in X)$ such that for every $x\in X$ the set $I_x=\{i\in I:x\in V_i\}$ contains at most two elements. Let $(\varphi_i:i\in I)$ be a partition of unity on $X$ such that $V_i=\varphi_i^{(-1)}((0,1])$ for every $i\in I$. For every  $i\in I$ we choose such  $x_i\in X$ that $V_i\subseteq U_{x_i}$ and let $y_i=y_{x_i}$. For every $x\in X$ let $g(x)=\sum\limits_{i\in I}\varphi_i(x)f(x,y_i)$. Clearly,  $g$ is continuous. Moreover, according to the choice of $U_x$ the following condition holds:

$(a)$\,\,\,\,$t+\frac{1}{2}\alpha(t)<y_i<t+\alpha(t)$ for every $t\in V_i$.

For every $x\in X$ we {pick} $i_x\in I_x$. Moreover, we take  $j_x\in I_x\setminus \{i_x\}$ if $|I_x|=2$, and $j_x=i_x$ if $|I_x|=1$. Let $\beta(x)=y_{i_x}-x$ and $\gamma(x)=y_{j_x}-x$. Since $x\in V_{i_x}$ and $x\in V_{j_x}$, by $(a)$ we have $\frac{1}{2}\alpha(x)<\beta(x) <\alpha(x))$ and $\frac{1}{2}\alpha(x)<\gamma(x) <\alpha(x))$.

Let $|I_x|=1$, i.e. $I_x=\{i\}$. Then $i_x=i$ and $g(x)=f(x,y_{i})=f(x,x+\beta(x))$. Now {assume} $|I_x|=2$. Then $I_x=\{i_x,j_x\}$ and $$p(g(x)-f(x,x+\beta(x)))=p\left(\sum\limits_{i\in I_x}\varphi_i(x)f(x,y_i)-f(x,y_{i_x})\right)=$$ $$=p(\varphi_{j_x}(x)f(x,y_{j_x})-\varphi_{j_x}(x)f(x,y_{i_x}))\leq $$
$$
\le p(f(x,y_{j_x})-f(x,y_{i_x}))=p(f(x,x+\beta(x))-f(x,x+\gamma(x))).$$
\hfill$\Box$

\begin{theorem}\label{th:3.2}
  Let $X=\mathbb R$, let $Z$ be a metric linear space with the metric,  generated by a pseudo-norm $p$,  and let $f:X^2\to Z$ be a mapping which is continuous with respect to the first variable and has a finite variation and is continuous with respect to the second variable at every point of the diagonal  $\Delta=\{(x,x):x\in X\}$. Then $g(x)=f(x,x)$ is an absolute Baire-one mapping.
\end{theorem}

{\bf Proof.} For every $n\in\mathbb N$ we apply Proposition~\ref{pr:3.1} to $f$ and to $\alpha_n:X\to(0,+\infty)$, where $\alpha_n(x)=\frac{1}{2^n}$. Then we obtain sequences $(g_n)^{\infty}_{n=1}$, $(\beta_n)^{\infty}_{n=1}$ and $(\gamma_n)^{\infty}_{n=1}$ of continuous mappings  $g_n:X\to Z$ and functions $\beta_n:X\to (\frac{1}{2^{n+1}},\frac{1}{2_n})$ and $\gamma_n:X\to (\frac{1}{2^{n+1}},\frac{1}{2_n})$ such that $$p(g_n(x)-f(x,x+\beta_n(x)))\leq p(f(x,x+\beta_n(x))-f(x,x+\gamma_n(x)))$$ for any $n\in\mathbb N$ and $x\in X$.

We show that the sequence $(g_n)^{\infty}_{n=1}$ is as required. Fix $x\in X$ and let $z_n=g_n(x)$, $u_n=f(x,x+\beta_n(x))$ and $v_n=f(x,x+\gamma_n(x))$ for every $n\in\mathbb N$. Since $\lim\limits_{n\to\infty}\beta_n(x)=\lim\limits_{n\to\infty}\gamma_n(x)=0$ and $f$ is continuous with respect to the second variable at $(x,x)$, $\lim\limits_{n\to\infty}u_n= \lim\limits_{n\to\infty}v_n=f(x,x)$ and $\lim\limits_{n\to\infty}p(u_n-v_n)=0$. Now, taking into account that  $p(z_n-u_n)\leq p(u_n-v_n)$ for every $n\in\mathbb N$, we obtain that  $\lim\limits_{n\to\infty}p(z_n-u_n)=0$ and $$\lim\limits_{n\to\infty}g_n(x)=\lim\limits_{n\to\infty}z_n=\lim\limits_{n\to\infty}(z_n-u_n)+\lim\limits_{n\to\infty}u_n=f(x,x)=g(x).$$ Remark that for every $n\in\mathbb N$ the points $x+\beta_n(x)$ and $x+\gamma_n(x)$ belong to the interval $I_n=(x+\frac{1}{2^{n+1}},x+\frac{1}{2_n})$ and  $I_n\cap I_m=\O$ for all distinct $n$ and $m$. Moreover, $\lim\limits_{n\to\infty}(x+\beta_n(x))=\lim\limits_{n\to\infty}(x+\gamma_n(x))=x$ and the mapping $f^x:X\to Z$, $f^x(t)=f(x,t)$, has a finite variation at  $x$. Hence, $$\sum\limits_{n=1}^{\infty}p(f^x(x+\beta_n(x))-f^x(x+\gamma_n(x)))=\sum\limits_{n=1}^{\infty}p(u_n-v_n)=C_1<\infty$$ and $$\sum\limits_{n=1}^{\infty}p(f^x(x+\beta_{n+1}(x))-f^x(x+\beta_n(x)))=\sum\limits_{n=1}^{\infty}p(u_{n+1}-u_n)=C_2<\infty.$$ Taking into account that $p(z_n-u_n)\leq p(u_n-v_n)$ for all $n\in\mathbb N$, we have $$\sum\limits_{n=1}^{\infty}p(g_{n+1}(x)-g_n(x))=\sum\limits_{n=1}^{\infty}p(z_{n+1}-z_n)$$
$$\leq \sum\limits_{n=1}^{\infty}\left(p(z_{n+1}-u_{n+1}) +  p(u_{n+1}-u_n)+ p(z_n -u_n)\right)\leq$$ $$\leq \sum\limits_{n=2}^{\infty}p(u_n-v_n)  + \sum\limits_{n=1}^{\infty}p(u_{n+1}-u_n) + \sum\limits_{n=1}^{\infty}p(u_n-v_n)\leq 2C_1+C_2<\infty.$$ Thus, $\sum\limits_{n=1}^{\infty}|g_{n+1}(x)-g_n(x)|_Z<\infty$ and $g$ is an absolute Baire-one mapping.
\hfill$\Box$

\begin{corollary}\label{cor:3.3}
 Let $X\subseteq\mathbb R$ be an interval, let $Z$ be a {metric linear space} with the metric, generated by a pseudo-norm, and let $f:X^2\to Z$ be a separately continuous mapping, which {has finite} variation with respect to the second variable at every point of the diagonal $\Delta=\{(x,x):x\in X\}$. Then the mapping $g:X\to Z$, $g(x)=f(x,x)$, is absolute Baire-one mapping.
\end{corollary}

{\bf Proof.} a). Let  $X=(a,b)$, where $a\in \mathbb R\cup\{-\infty\}$ and $b\in \mathbb R\cup\{+\infty\}$. Consider a homeomorphism $\varphi:(a,b)\to\mathbb R$ such that $\varphi$ and $\varphi^{-1}$ are locally Lipschitz.
Then the mapping $\tilde{f}:\mathbb R\to Z$, $\tilde{f}(t)=f(\varphi^{-1}(t))$, satisfies the conditions of Theorem~\ref{th:3.2}. Therefore, there exists a sequence $(\tilde{g}_n)_{n=1}^{\infty}$ of continuous functions $\tilde{g}_n:\mathbb R\to Z$ such that $\lim\limits_{n\to\infty}\tilde{g}_n(t)=\tilde{f}(t,t)$ and $\sum\limits_{n=1}^{\infty}|\tilde{g}_{n+1}(t)-\tilde{g}_n(t)|_Z<\infty$ for every $t\in \mathbb R$. Then the sequence $(g_n)_{n=1}^{\infty}$ of continuous functions $g_n:X\to Z$, $g_n(x)=\tilde{g}_n(\varphi(x))$, {pointwise converges to} $g$ on $X$ and $\sum\limits_{n=1}^{\infty}|g_{n+1}(x)-g_n(x)|_Z<\infty$ for every $x\in X$.

b). Now let $X=[a,b)$. According to a) there is a sequence $(h_n)_{n=1}^{\infty}$ of continuous functions $h_n:(a,b)\to Z$ such that $\lim\limits_{n\to\infty}h_n(x)=g(x)$ and $\sum\limits_{n=1}^{\infty}|h_{n+1}(x)-h_n(x)|<\infty$ for every $x\in (a,b)$. We choose a convergent to $a$ monotone sequence $(a_n)_{n=1}^{\infty}$ of reals $a_n\in(a,b)$ and for every $n\in\mathbb N$ let $g_n:X\to Z$ be a continuous mapping such that  $g_n(a)=g(a)$ and $g_n(x)=h_n(x)$ for every $x\in[a_n,b)$. Then $(g_n)_{n=1}^{\infty}$ {pointwise converges to} $g$ on $X$ and  $\sum\limits_{n=1}^{\infty}|g_{n+1}(x)-g_n(x)|<\infty$ for every $x\in X$.

c). In the case when $X=(a,b]$ or $X=[a,b]$ we argue analogously {as we did to prove} b).\hfill$\Box$

\section{The characterization of diagonals of separately absolutely continuous functions}

\begin{proposition}\label{pr:5.1}
 Let $X$ be a metric space, let $Z$ be a normed space, let $f:X\to Z$ be a continuous mapping and $\varepsilon>0$. Then there exists a locally Lipschitz mapping $g:X\to Z$ such that $\|f(x)-g(x)\|\leq \varepsilon$ for every $x\in X$.
\end{proposition}

{\bf Proof.} Since $f$ is continuous and $X$ is a paracompact by \cite[Theorem 5.1.3]{Eng}, there exists a locally finite cover $\mathcal U$ of $X$ by open non-empty sets $U$ with ${\rm diam}\, f(U)<\varepsilon$. For every $U\in\mathcal U$ and $x\in X$ by $\psi_U(x)$ we denote the distance from $x$ to $X\setminus U$. Let $\varphi_U=\frac{\psi_U}{\sum_{V\in\mathcal U}\psi_V}$ for every $U\in\mathcal U$. Note that $(\varphi_U:U\in\mathcal U)$ is a partition of unity on $X$, all the functions  $\varphi_U:X\to[0,1]$ are locally Lipschitz and ${\rm supp}\,\varphi_U=U$. Moreover, for every  $U\in\mathcal U$ we take $z_U\in f(U)$. Consider the function $g:X\to Z$, $g(x)=\sum\limits_{U\in\mathcal\, U}\varphi_U(x)z_U$. Is is easy to see that $\|f(x)-g(x)\|\leq \varepsilon$ for every $x\in X$. Since $\mathcal U$ is locally finite and $\varphi_U$ are locally Lipschitz, $g$ is locally Lipschitz.
\hfill$\Box$

The following example shows that the {analogue} of Proposition~\ref{pr:5.1} for locally convex $F$-spaces $Z$ is not valid.

\begin{example}\label{ex:5.2}
 Let $0<p<1$, let $X=\mathbb R$ be a space with Euclid metric, let $Z= \mathbb R$ be a space with the metric $|x-y|_Z=|x-y|^p$,  $f(x)=x$ for every $x\in X$ and $\varepsilon>0$. Since every locally Lipschitz mapping  $g:X\to Z$ is a constant,  $\sup \{ |f(x)-g(x)|_p:x\in X\}>\varepsilon$ for any locally Lipschitz mapping  $g:X\to Z$.
\end{example}

\begin{proposition}\label{pr:5.3}
  Let $g:[a,b]\to\mathbb R$ and $h:[a,b]\to\mathbb R$ be Lipschitz functions with a constant $K$, let $\varphi :[a,b]\to[0,1]$ be a Lipschitz function with a constant $L$ and $|g(x)-h(x)|\leq M$ for every $x\in [a,b]$. Then the function $$f(x)=\varphi(x)g(x)+(1-\varphi(x))h(x)$$ is Lipschitz with the constant $C=K+LM$ on $[a,b]$.
\end{proposition}

 {\bf Proof.} Let $x_1,x_2\in [a,b]$. {Then we} have $$|f(x_2)-f(x_1)|= |(\varphi(x_2)-\varphi(x_1))(g(x_2)-h(x_2))+\varphi(x_1)(g(x_2)-g(x_1))+$$ $$+(1-\varphi(x_1))(h(x_2)-h(x_1))|\leq L M |x_2-x_1|+ \varphi(x_1)K|x_2-x_1|+$$ $$+(1-\varphi(x_1))K|x_2-x_1|= (K+LM)|x_2-x_1|.$$\hfill$\Box$

\begin{proposition}\label{pr:5.4}
  Let $Z$ be a metric space, $f:[a,b]\to\mathbb Z$ a continuous mapping and let $(a_n)^{\infty}_{n=1}$ be a strictly increasing sequence such that $a_1=a$, $\lim\limits_{n\to\infty}a_n=b$, let for every  $n\in\mathbb N$ the mapping $f$ is Lipschitz with a constant  $C_n$ on $[a_n,a_{n+1}]$, {and} moreover, the series $\sum\limits_{n=1}^{\infty}C_n(a_{n+1}-a_n)$ is convergent. Then $f$ is absolutely continuous on $[a,b]$.
\end{proposition}

{\bf Proof.}  Fix  $\varepsilon>0$ and choose $m\in\mathbb N$ such that $\sum\limits_{n>m}C_n(a_{n+1}-a_n)<\frac{\varepsilon}{2}$. Set $\displaystyle \delta=\frac{\varepsilon}{2\max\{C_n:1\leq n\leq m\}}$. Let $$a\leq u_1<v_1\leq u_2<v_2\leq\dots \leq u_i<v_i\leq a_m\leq
$$
 $$
 \le u_{i+1}<v_{i+1}\leq\dots\leq u_{i+j}< v_{i+j}\leq b$$ with $\sum\limits_{k=1}^{i+j}(v_k-u_k)<\delta$. Then $$\sum\limits_{k=1}^{i+j}|f(v_k)-f(u_k)|_Z=\sum\limits_{k=1}^{i}|f(v_k)-f(u_k)|_Z+\sum\limits_{k=i+1}^{i+j}|f(v_k)-f(u_k)|_Z\leq$$ $$\leq\max\{C_n:1\leq n\leq m\}\cdot\sum\limits_{k=1}^{i}(v_k-u_k)+\sum\limits_{n>m}C_n(a_{n+1}-a_n)<\frac{\varepsilon}{2}+\frac{\varepsilon}{2}=\varepsilon.$$

 Now let $$a\leq u_1<v_1\leq \dots \leq u_i<v_i\leq b$$ with $\sum\limits_{k=1}^{i}(v_k-u_k)<\delta$ and $k\leq i$ such that $a_m\in(u_k,v_k)$. Using the above-obtained
 estimation we have
 $$\sum\limits_{j=1}^{i}|f(v_j)-f(u_j)|_Z=\sum\limits_{j=1}^{k-1}|f(v_j)-f(u_j)|_Z+ |f(a_m)-f(u_k)|+$$
 $$|f(v_k)-f(a_m)|+\sum\limits_{j=k+1}^{i}|f(v_j)-f(u_j)|_Z\leq\varepsilon.$$ \hfill$\Box$

\begin{theorem}\label{th:5.5}
   Let $X\subseteq \mathbb R$ be an interval, let $Z$ be a normed space and let  $g:X\to Z$ be an absolute Baire-one mapping. Then there exists a separately absolutely continuous mapping  $f:X^2\to Z$ such that  $g(x)=f(x,x)$ for every $x\in X$.
\end{theorem}

{\bf Proof.} Let $(I_n)_{n=1}^{\infty}$ be an increasing sequence of segments $I_n=[a_n,b_n]$ such that $X=\bigcup\limits_{n=1}^{\infty}I_n$. Let $g_0(x)=0$ for every $x\in X$. Using the definition of absolute Baire-one class and Proposition~\ref{pr:5.1}, we choose a sequence $(\tilde{g}_n)_{n=1}^{\infty}$ of locally Lipschitz mappings  $\tilde{g}_n:X\to Z$ such that $\lim\limits_{n\to\infty}\tilde{g}_n(x)=g(x)$ and $\sum\limits_{n=1}^{\infty}\|\tilde{g}_{n+1}(x)-\tilde{g}_n(x)\|<\infty$ for every $x\in X$. Now for arbitrary  $n\in\mathbb N$ and $x\in X$ let
$$
 g_n(x)=\left\{\begin{array}{lll}
                         \tilde{g}_n(x), & x\in I_n\\
                         \tilde{g}_n(a_n), & x\in X\cap(-\infty,a_n)\\
                         \tilde{g}_n(b_n), & x\in X\cap(b_n,+\infty).
                       \end{array}
 \right.
 $$
Then we obtain the sequence $(g_n)_{n=1}^{\infty}$ of Lipschitz mappings $g_n:X\to Z$ such that $\lim\limits_{n\to\infty}g_n(x)=g(x)$ and $\sum\limits_{n=1}^{\infty}\|{g}_{n+1}(x)-{g}_n(x)\|<\infty$ for every  $x\in X$.

Let $(K_n)^{\infty}_{n=1}$ be a sequence of reals $K_n>0$ such that the functions $g_{n-1}$ and $g_{n}$ are Lipschitz with the constant $K_n$. We choose a strictly decreasing sequence $(\delta_n)^{\infty}_{n=1}$ of reals $\delta_n>0$ such that $\lim\limits_{n\to\infty}\delta_n=0$ and $\sum\limits_{n=1}^{\infty}K_n\delta_n<\infty$.

For every $n\in\mathbb N$ let $F_n=\{(x,y)\in X^2:|x-y|\leq \delta_n\}$, $G_n=\{(x,y)\in X^2:|x-y|< \delta_n\}$. {We consider} the function $\varphi_n:\mathbb R\to[0,1]$ {defined by}
$$
 \varphi_n(t)=\left\{\begin{array}{lll}
                         1, & |t|>\delta_n\\
                         \frac{t-\delta_{n+1}}{\delta_n-\delta_{n+1}}, & \delta_{n+1}\leq |t|\leq \delta_n\\
                         0, & |t|<\delta_{n+1}.
                       \end{array}
 \right.
 $$
Remark that every function $\varphi_n$ is Lipschitz with the constant $L_n=\frac{1}{\delta_n-\delta_{n+1}}$.

We consider the following function $f:X^2\to Z$:
$$
 f(x,y)=\left\{\begin{array}{lll}
                         g_0(x), & |x-y|>\delta_1\\
                         \varphi_n(x-y)g_{n-1}(x)+(1-\varphi_n(x-y))g_n(x), & (x,y)\in F_n\setminus F_{n+1}\\
                         g(x), & x=y\in X.
                       \end{array}
 \right.
 $$
Let us show that $f$ {has the desired properties}. Since for every $n\in\mathbb N$ with $|x-y|=\delta_{n+1}$ we have $\varphi_n(x-y)=0$ and $\varphi_{n+1}(x-y)=1$, $$f(x,y)=g_n(x)=\varphi_{n}(x-y)g_{n-1}(x)+(1-\varphi_{n}(x-y))g_{n}(x).$$ Therefore, $$f(x,y)= \varphi_n(x-y)g_{n-1}(x)+(1-\varphi_n(x-y))g_n(x)$$ for every $(x,y)\in F_n\setminus G_{n+1}$. The continuity of $g_n$ and $\varphi_n$ implies that $f$ is jointly continuous on $X^2\setminus \Delta$, where $\Delta=\{(x,x):x\in X\}$. Moreover, the equality $\lim\limits_{n\to\infty}g_n(x)=g(x)$ fir every $x\in X$ implies that $f$ is continuous with respect to the first and to the second variable at all points of $\Delta$. Thus, $f$ is separately continuous. Moreover, since all of $g_n$ and all of $\varphi_n$ are Lipschitz,  $f$ is Lipschitz on  $X^2\setminus G_n$ for every $n\in\mathbb N$.

{We fix} $x_0\in X$  and show that the mapping $f_{x_0}:X\to Z$, $f_{x_0}(x)=f(x,x_0)$, is absolutely continuous. It is sufficient to prove that the mapping $f_{x_0}$ is absolutely continuous on the intervals  $X_1=X\cap(-\infty,x_0]$ and $X_2=X\cap [x_0,+\infty)$. Let us consider the case when the interval $X_1$ is non-degenerated, i.e. $X_1\cap(-\infty,x_0)\ne\emptyset$. Since $X=\bigcup\limits_{n=1}^{\infty}I_n$ and $(I_n)^{\infty}_{n=1}$ increases and $\lim\limits_{n\to\infty}\delta_n=0$, there exists a number $m\in\mathbb N$ such that  $[x_0-\delta_m,x_0]\subseteq I_m$. Observe that $f_{x_0}(x)=0$ for every $x\in X\cap (-\infty,x_0-\delta_1]$ and
$f_{x_0}$ is Lipschitz on $X_1\cap [x_0-\delta_1,x_0-\delta_m]$. Therefore, the absolutely continuity of $f_{x_0}$ of $X_1$  is equivalent to the absolutely continuity of $f_{x_0}$ on  $[x_0-\delta_m, x_0]$.

Fix $n\geq m$. Notice that $g_{n-1}$ and $g_n$ are Lipschitz with the constant $K_n$, and $\varphi_n$ is Lipschitz with the constant  $L_n$. Moreover, for every $x\in [x_0-\delta_n,x_0-\delta_{n+1}]$ we have
$$
\|g_{n}(x)-g_{n-1}(x)\|\leq \|g_{n}(x_0)-g_{n-1}(x_0)\|+\|g_{n}(x)-g_{n}(x_0)\|+$$
$$
+\|g_{n-1}(x)-g_{n-1}(x_0)\|\leq\|g_{n}(x_0)-g_{n-1}(x_0)\|+2K_n\delta_n=M_n.
$$
Hence,  $f_{x_0}$ is Lipschitz on $[x_0-\delta_n,x_0-\delta_{n+1}]$ with the constant $$C_n=K_n+L_nM_n=K_n+\frac{1}{\delta_n-\delta_{n+1}}(\|g_{n}(x_0)-g_{n-1}(x_0)\|+2K_n\delta_n).$$ by Proposition~\ref{pr:5.3}. Now we have $$\sum\limits_{n\geq m}C_n(\delta_n-\delta_{n+1})\leq 3\sum\limits_{n\geq m}K_n\delta_n+ \sum\limits_{n\geq m}\|g_n(x_0)-g_{n-1}(x_0)\|<\infty\,.$$ By Proposition~\ref{pr:5.4}, $f_{x_0}$ is absolutely continuous in $[x_0-\delta_m, x_0]$, consequently, on $X_1$.

The absolutely continuity of $f_{x_0}$ on $X_2$ and the absolutely continuity of $f$ with respect to the second variable are proved similarly.
\hfill$\Box$

 Theorems~\ref{th:3.2} and~\ref{th:5.5} imply the following result.

\begin{theorem}\label{th:5.6}
 Let $X\subseteq\mathbb R$ be an interval, let $Z$ be a normed space and $g:X\to Z$. Then the following conditions are equivalent:

 $(i)$\,\,\,there exists a separately absolute continuous mapping $f:X^2\to Z$ with the diagonal  $g$;

 $(ii)$\,\,\,there exists a mapping $f:X^2\to Z$ with the diagonal $g$ which is continuous with respect to the first variable and has a finite variation and is continuous with respect to the second variable at every point of the diagonal \mbox{$\Delta=\{(x,x):x\in X\}$};

$(iii)$\,\,\,$g$ is an absolute Baire-one mapping.
\end{theorem}

\section{Open problems}

\begin{question}\label{q:6.2}
 Let $X=[0,1]$, let $Z$ be a metric space and let $f:X^2\to Z$ be a separately absolute continuous mapping.  Is  $g(x)=f(x,x)$ an absolute Baire-one mapping?
\end{question}

\begin{question}\label{q:6.3}
 Let $X=[0,1]$, let $Z$ be a metric topological vector space and let  $g:X\to Z$ be an absolute Baire-one mapping. Does there {exist a separately} absolute continuous mapping $f:X^2\to Z$ such that $g(x)=f(x,x)$ for every $x\in X$?
\end{question}

\end{document}